\documentclass{birkjour}

 \newtheorem{thm}{Theorem}[section]
 \newtheorem{cor}[thm]{Corollary}
 
 \newtheorem{prop}[thm]{Proposition}
 \theoremstyle{definition}
 \newtheorem{defn}[thm]{Definition}
 \theoremstyle{remark}

 \numberwithin{equation}{section}
 \usepackage{color}
\DeclareMathOperator{\Imagine}{Im}
 \usepackage{hyperref}
 
\begin{document}

%
%
%
%

\title[Suita conjecture for a punctured torus]
 {Suita conjecture for a punctured torus}

\author[R. X. Dong]{Robert Xin DONG}
\address{Instytut Matematyki, Uniwersytet Jagiello\'nski, \L{}ojasiewicza 6, 30-348 Krak\'ow, Poland}
\email{1987xindong@tongji.edu.cn}
\dedicatory{To my family}

\begin{abstract}
For a once-punctured complex torus, we compare the Bergman kernel and the fundamental metric, by constructing explicitly the Evans-Selberg potential and discussing its asymptotic behaviors. This work aims to generalize the Suita type results to potential-theoretically parabolic Riemann surfaces.
\end{abstract}

\maketitle

\section {Introduction}
The Suita conjecture \cite{Su} asks about the precise relations between the Bergman kernel and the logarithmic capacity. For potential-theoretically hyperbolic Riemann surfaces, it was conjectured that the Gaussian curvature of the Suita metric (induced from the logarithmic capacity) is bounded from above by $-4$. The relations between the Suita conjecture and the $L^2$ extension theorem were first observed in \cite{Oh95} and later contributed by several mathematicians. The Suita conjecture proved to be true for the hyperbolic case (see \cite{Bl12, G-Z15, BL, O}), and it might be interesting to generalize similar results to \verb+non+-hyperbolic cases. In this article, for a once-punctured complex torus $X_{\tau, u}:=X_\tau \backslash \{u\}$, which is a typical potential-theoretically parabolic Riemann surface, we construct a so-called Evans-Selberg potential and further derive a so-called fundamental metric.

\begin{thm} \label{punctured} On $X_{\tau, u}$, let $K_{\tau, u}$ and $c_{\tau, u}$ be the Bergman kernel and the fundamental metric, respectively. In the local coordinate $w$, write $K_{\tau, u}=k_{\tau, u}(w)|dw|^2$ and $c_{\tau, u}=c_{\tau, u}(w)|dw|^2$. Then, as $w\to u$, $$\frac{\pi k_{\tau, u}(w)}{c_{\tau, u}^2(w)} \sim  \frac{ \pi\cdot  |w-u|^2 }{2 \cdot \Imagine \tau }\to 0^+.$$ \end{thm} 

We also study the degenerate case when a once-punctured complex torus becomes a singular curve and obtain the following result.

\begin{thm} \label{degenerate} Under the same assumptions as in Theorem \ref{punctured}, as $\Imagine \tau \to +\infty,$ it follows that $$\frac{\pi K_{\tau, u}(w)}{c^2_{\tau, u}(w)} \to 0^+.$$
\end{thm}

As we can see, either Theorem \ref{punctured} or Theorem \ref{degenerate} will imply that the Gaussian curvature of the fundamental metric on $X_{\tau, u}$ can be arbitrarily close to $0^-$, which is different from the hyperbolic case. 

\begin{cor} The Gaussian curvature of the fundamental metric on a once-punctured complex torus cannot be bounded from above by a negative constant.
\end{cor}

\section{Preliminaries} On a potential-theoretically parabolic Riemann surface $R$, there exist a so-called Evans-Selberg potential, which is a counterpart of the Green function for the hyperbolic case. Let us recall the definition of an Evans-Selberg potential and the so-called fundamental metric (cf. \cite[p.351]{SN}, \cite [p.114]{SNo}, \cite{McV}). 

\begin {defn}\label{def-Evans} On an open Riemann surface $\Sigma$, an Evans-Selberg potential $E_{q}(p)$ with a pole $q\in \Sigma$ is a real-valued function such that:

$(i)$ For all $p\in \Sigma\setminus \{q\}$, $E_{q}(p)$ is harmonic with respect to $p$,  

$(ii)$ ${ E_{q}(p)} \to +\infty$, as $p\to a_{\infty}$ (the Alexandroff ideal boundary point),

$(iii)$ $E_{q}(p)-\log|\varphi(p)-\varphi(q)|$ is bounded near $q$, with $\varphi$ being the local coordinate.
\end {defn} 

\begin{defn} On a potential-theoretically parabolic Riemann surface $\Sigma$, the fundamental metric under the local coordinate $z=\varphi(p)$ is defined as \begin{equation*}\label{fundamental metric}
c(z)|dz|^2:=\exp \lim_{q \to p}\left(E_q(p)-\log|\varphi(p)-\varphi(q)|\right)|dz|^2.
\end{equation*} \end{defn}

The fundamental metric is a non-compact counterpart of the Arakelov metric, and it coincides at the hyperbolic case with the Suita metric. The Gaussian curvature form of the fundamental metric is \begin{equation} \label{varolin} -4\frac{\partial^2}{\partial z \partial \bar z}\log c(z)= -4 \pi k(z),\end {equation} where $k(z)$ $(\geq 0)$ is the coefficient of the Bergman kernel $(1, 1)$-form in the local coordinate $z$. For a compact complex torus $X_\tau:=\mathbb C/\left(\mathbb Z+\tau\mathbb Z\right)$, where $\tau\in \mathbb C$ and $\Imagine \tau>0$, its Bergman kernel by definition is $K_{\tau}(z)=(\Imagine \tau)^{-1} dz\wedge d\bar z,$ where $z$ is the local coordinate induced from the complex plane $\mathbb C$.

\begin{defn} A once-punctured complex torus $X_{\tau, u}:=X_\tau \backslash \{u\}$ is an open Riemann surface obtained by removing one single point $u$ from a compact complex torus $X_\tau$. \end{defn} 

 \begin{prop} There exists no non-constant subharmonic function which is bounded from above on $X_{\tau, u}$.\end{prop}

On the one hand, the parabolicity of $X_\tau$ follows straightforwardly from Removable Singularity Theorem, which also applies to finitely-punctured Riemann surfaces. On the other hand, there exists a proof (by using Maximum Principle and finding a harmonic function, as we will see below) that works for an open Riemann surface $X$ of infinitely many genus. Equivalently, we could say that $X$ admits an Evans-Selberg potential.

\begin{prop} \label{infinite} There exists no non-constant subharmonic function which is bounded from above on the algebraic curve $$X:=\left\{(y,x)\in \mathbb C^2\, \left |\, y^2=x \prod_{n=1}^{\infty}\left(1- {x^2}/{n^2}\right)\right.\right\}.$$
\end{prop}

\begin{proof} [Sketch of proof] The idea is by contradiction and suppose there exists a subharmonic function $u$ which is bounded from above on $X$. On $X$, we define \begin{displaymath}
f(y, x):=\left\{ \begin{array}{ll}0, &|x|\leq1\\
\log |x|, &|x|>1,\\
\end{array} \right.
\end{displaymath} and consider $v:=u-\epsilon f$. Taking a limit as $\epsilon$ tends to $0$, one will prove that the maximum is attainable at the ``lift" of the unit circle (away from the boundary). By Maximum Principle, one concludes that $u$ is constant. \end{proof}

Explicit formulas of Evans-Selberg potentials and fundamental metrics on planar domains are provided in \cite{D17}.

\section{A Compact Torus} The Arakelov-Green function on a complex torus $X_\tau$ with a pole $w$ satisfies \begin {equation} \label{AG} \frac{\partial ^2 g_w(z)}{\partial z \bar \partial z}=\frac{\pi}{2} \left( \delta(z-w)-\frac{1}{ \Imagine \tau}\right),\end {equation} and can be expressed via the theta function as \begin {equation} \label{green} g_w(z)=\log \left | \frac { \theta_1 (z-w; q)} {\eta (\tau)}\right|-\frac{\pi \cdot (\Imagine (z-w))^2}{\Imagine \tau}.\end {equation} 

Here $\eta (\tau)=q^{\frac{1}{12}} \cdot \prod_{m=1}^{\infty}(1-q^{2m})$ and $$\theta_{1}(z;q):=2q^{1/4}\sin (\pi z)\prod_{n=1}^\infty (1-q^{2n})(1-2\cos(2\pi z)q^{2n}+q^{4n}),$$ for $q=\exp(\pi i \tau)$. In this case, it is possible to compare the Bergman kernel and the Arakelov metric. For $X_\tau$, the author in \cite{D14} numercially obtained a sharp negative upper bound for the Arakelov metric by computing via elliptic functions\footnote{The author apologizes for several mistakes contained in \cite{D14}.}. Alternatively, the Gaussian curvature of the Arakelov metric can be computed by \textit{Mathematica} (Version 10.3). Recall that for the hyperbolic case, the Green function is always less than $0$ in the interior. However, the supremum of an Arakelov-Green function on a torus can be positive at some points, illustrated in Figure 1 (plotted by \textit{Mathematica}). 

 \begin{figure}[htp]
  \centering
  \includegraphics[width=0.75\textwidth]{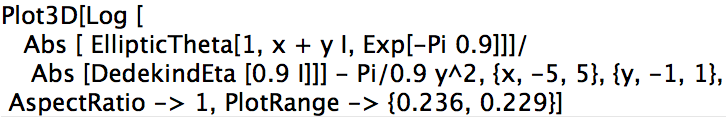}
\end{figure}

 \begin{figure}[htp]
  \centering
  \includegraphics[width=0.75\textwidth]{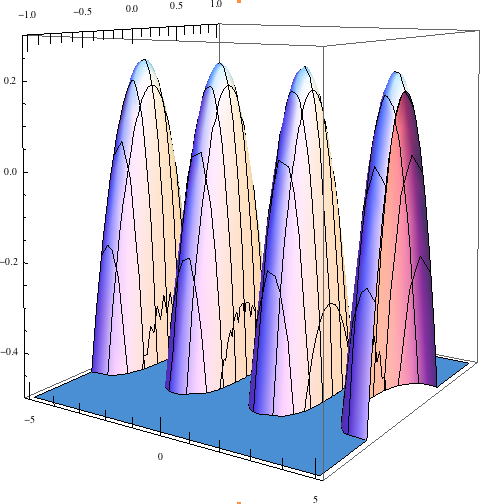}
   \caption{An Arakelov-Green function on a torus}
\end{figure}

\section{A Once-Punctured Torus} 

\begin{thm} \label{Evans} There exists an Evans-Selberg potential on $X_{\tau, u}$ with a pole $w$ given by $$E^{\tau, u}_w(z)=\log \left | \frac { \theta_1 (z-w; q)} {\theta_1 (z-u; q)}\right|,$$ for $z\in X_{\tau, u}\setminus \{w\}.$\end{thm}

\begin{proof} [Proof of Theorem \ref{Evans}] We see that the two terms on the right-hand side of \eqref{green} are responsible for the two terms on the right-hand side of \eqref{AG}, respectively. Keeping this in mind, we can construct the Evans-Selberg potential by attaching physics meanings. We regard the potential as an electric flux generated at the pole $w$ and terminates at the boundary point $u$ (see \cite{Oo10} for detailed physics explanations). So, the Evans-Selberg potential $E^{\tau, u}_w(z)$ with a pole $w$ satisfies $$\frac{\partial ^2 E^{\tau, u}_w(z)}{\partial z \bar \partial z}=\frac{\pi}{2} \left( \delta(z-w)- \delta(z-u)\right),$$ and can be expressed via the theta function as \begin{align*} \hspace{-0.38cm} E^{\tau, u}_w(z)&=\log \left | \frac { \theta_1 (z-w; q)} {\eta (\tau)}\right|- \log \left | \frac { \theta_1 (z-u; q)} {\eta (\tau)}\right|=\log \left | \frac { \theta_1 (z-w; q)} {\theta_1 (z-u; q)}\right|\\
&=\log \left| \frac{\sin(\pi(z-w))\cdot \prod_{m=1}^{\infty}(1-2\cos(2\pi (z-w))\cdot q^{2m}+q^{4m})}{\sin(\pi(z-u))\cdot \prod_{m=1}^{\infty}(1-2\cos(2\pi (z-u))\cdot q^{2m}+q^{4m})}\right|.\end{align*} \end{proof}

\begin{cor} \label{fund} There exists a fundamental metric $c_{\tau, u}$ on $X_{\tau, u}$ in the local coordinate $w$ given by $$c_{\tau, u}(w)|dw|^2=\frac{ 2\pi \cdot |\eta(\tau)|^3 }{   |\theta_1 (w-u; q)|}|dw|^2.$$ \end{cor}

\noindent \textit{Proof.} This can be verified by definition, since \begin{align*}
\hspace{0cm}   c_{\tau, u}(w)&=\exp \lim_{z \to w}\left(\log \left | \frac { \theta_1 (z-w; q)} {\theta_1 (z-u; q)}\right|-\log |z-w|\right)\\
& =  \left| \frac{ \pi \cdot \prod_{m=1}^{\infty}(1-q^{2m})^2 }{\sin(\pi(w-u))\cdot \prod_{m=1}^{\infty}(1-2\cos(2\pi (w-u))\cdot q^{2m}+q^{4m})}\right|\\
&=\left| \frac{ \pi \cdot \eta(\tau)^2  }{q^{\frac{1}{6}}\cdot\sin(\pi(w-u))\cdot \prod_{m=1}^{\infty}(1-2\cos(2\pi (w-u))\cdot q^{2m}+q^{4m})}\right|\\
&=\left| \frac{ \pi \cdot \eta(\tau)^2\cdot 2q^{\frac{1}{4}} \cdot\prod_{m=1}^{\infty} (1-q^{2m}) }{q^{\frac{1}{6}}\cdot \theta_1 (w-u; q)}\right|\\
&=\left| \frac{ \pi \cdot \eta(\tau)^2\cdot 2q^{\frac{1}{4}} \cdot \frac{\eta(\tau)}{q^{\frac{1}{12}}} }{q^{\frac{1}{6}}\cdot \theta_1 (w-u; q)}\right|= \left| \frac{ 2\pi \cdot \eta(\tau)^3 }{   \theta_1 (w-u; q)}\right|.\tag*{\qed} \end{align*} 

By the second equality above, $c_{\tau, u}$ has the following asymptotic behavior, which will yield Theorem 1.2 for any fixed $\tau$.

\begin{cor} Under the same assumptions as in Corollary 4.2, as $w\to u$, it follows that $$ c_{\tau, u}(w)\sim \frac{ 1}{ |w-u| }\to +\infty.$$ \end{cor}  

\section{The Degenerate Case} By studying the asymptotic behaviors of the fundamental metric under degeneration with respect to the complex structure, we will prove Theorem \ref{degenerate}. Relating Theorem \ref{degenerate} with \eqref{varolin}, we further get Corollary 1.3.

\begin{proof} [Proof of Theorem \ref{degenerate}] Let $\Imagine \tau \to +\infty,$ and $q\equiv\exp(\pi i \tau)$ will tend to $0$. Then, it holds that \begin{align*}
c_{\tau, u}(w) & \to \left| \frac{ \pi \cdot \prod_{m=1}^{\infty}(1-0^{2m})^2 }{\sin(\pi(w-u))\cdot \prod_{m=1}^{\infty}(1-2\cos(2\pi (w-u))\cdot 0^{2m}+0^{4m})}\right| \\
& \to   \frac{ \pi  }{|\sin(\pi(w-u)) |} .  
\end{align*} 

Therefore, it follows that $$\frac{\pi K_{\tau, u}(w)}{c_{\tau, u}^2(w)} \to  \frac{  |\sin(\pi(w-u))|^2 }{2  \cdot \Imagine \tau  \cdot \pi } \to 0^+,$$ since the denominator is uniformly bounded by $1$ for any fixed $w$. \end{proof}

On the one hand, at the degenerate case of potential-theoretically hyperbolic Riemann surfaces, we are not sure  whether Gaussian curvatures of the Suita metrics are still bounded from above by $-4$. On the other hand, for a compact complex torus, the Gaussian curvature of the Arakelov metric is always $0$ by the genus reason, although our earlier result in \cite{D14} shows that as $\Imagine \tau \to +\infty,$ $$\frac{\pi K_{\tau }(w)}{c_{\tau}^2(w)} \to +\infty.$$\\

\paragraph{Acknowledgment} The author expresses his gratitude to Prof. T. Ohsawa for his guidance and to Prof. H. Umemura for the communications on elliptic functions. He also thanks H. Fujino and X. Liu for the discussions.

This work is supported by the Ideas Plus grant 0001/ID3/2014/63 of the Polish Ministry of Science and Higher Education, KAKENHI and the Grant-in-Aid for JSPS Fellows (No. 15J05093).

\end{document}